\documentclass[11pt]{article}
\usepackage{amsgen, amsmath, amsfonts, amsthm, amssymb, enumerate,amscd}
\usepackage[all]{xy}

\numberwithin{equation}{section}

\newtheorem{defn}{Definition}[section]

\newtheorem{lem}[defn]{Lemma}
\newtheorem{thm}[defn]{Theorem}
     
\newtheorem{rem}[defn]{Remark}
\newtheorem{conjecture}[defn]{Conjecture}

\title{An analogue of cyclotomic units for products of elliptic curves}  
\author{Srinath Baba and Ramesh Sreekantan\\ McGill University and
  Tata Institute of Fundamental Research}

\newcommand {\ZZ}{{\Bbb Z}}

\newcommand {\X}{{X}}

\newcommand {\E}{{E}}

\newcommand {\Z}{{\mathcal Z}}
\newcommand {\C}{{\mathcal C}} 
\newcommand {\Q}{{\Bbb Q}}
\newcommand {\OO}{{\mathcal O}}

\newcommand {\M}{{\mathcal M}}

\newcommand {\D}{{\mathcal D}}
\newcommand {\R}{{\Bbb R}}
\newcommand {\CC}{{\mathcal S}}

\begin{document}

\maketitle

\begin{abstract}
  We construct  certain  elements in the   integral motivic cohomology
  group $H^3_{\M}(E  \times  E',\Q(2))_{\ZZ}$, where $E$ and  $E'$ are
  elliptic curves over  $\Q$. When $E$ is  not isogenous to $E'$ these
  elements are analogous to `cyclotomic    units' in real    quadratic
  fields  as they come  from modular  parametrisations of the elliptic
  curves. We  then find an  analogue of the  class  number formula for
  real quadratic fields. Finally we use  the Beilinson conjectures for
  $E \times E'$ to deduce them for products  of $n$ elliptic curves. A
  certain amount of this paper is expository in nature.
\end{abstract}

\section {Introduction}

\subsection{Beilinson's Conjecture}

Let $\X$  be a variety  over a number field $F$.   The well known {\bf
  Tate Conjecture} relates the order of the pole at a certain point of
a cohomological $L$-function  of $\X$ as a variety  over $K$, a number
field containing $F$,   to  the rank of   a  group  of  cycles  modulo
homological equivalence defined over $K$.   However, Tate did not give
any interpretation of the residue at that point.

Beilinson  \cite{be} made  the  following generalisation  of the  Tate
conjecture  (For simplicity we will state it over $\Q$, as that is the
case we will be concerned with):

\begin{conjecture}[Beilinson]

  Let $\X$ be a smooth projective variety  defined over $\Q$.  Let $i$
  be an even integer  and  $m=\frac{i}{2}$.  Let $L(\X,s)$ denote  the
  $L$-function defined by $H^i(\X)$.   Let $B^{m}(\X)_{\Q}$ denote the
  $\Q$-vector space generated  by the codimensional $m$ cycles  modulo
  homological equivalence. Then:

\begin{itemize}
  
\item (i)           $\widetilde{r}_{\D}(B^m(\X)_{\Q}               \oplus
  H^{i+1}_{\M}(\X_,\Q(m+1))_{\ZZ})$ induces a $\Q$ structure on

$H^{i+1}_{\D}(\X_{/ \R},\R(m+1))$

\item (ii) $ord_{s=m} L(\X,s)=dim_{\Q} H^{i+1}_{{\mathcal M}}(\X,\Q(m+1))_{\ZZ}$

\item (iii) $ord_{s=m+1} L(\X,s)= - dim_{\Q}(B^m(\X)_{\Q})$ {\bf (Tate)} 

\item (iv) $ L^*(\X,s)_{s=m} \sim_{\Q^*} c_{\X}(m)$

\end{itemize}

\end{conjecture}

Here $\widetilde{r}_{\D}$ is a     certain `thickened' regulator  map    which
generalises   the     usual  regulator  map      for  number   fields,
$H^{i+1}_{\M}(\X,\Q(m+1))_{\ZZ}$ is the `integral' motivic cohomology,
which  is  the  motivic cohomology of  a   regular proper model  if it
exists, otherwise there is  an unconditional definition due  to Scholl
\cite{sc}  ,   and $H^{i+1}_{\D}(\X_{/   \R},\R(m+1))$  is  the `Real'
Deligne cohomology  which    is a real   vector   space  of  dimension
$ord_{s=m}     L(\X,s)    -   ord_{s=m+1}     L(\X,s)$     over  $\R$. 
$L^*(\X,s)_{s=m}$  is the first non-zero  term in the Taylor expansion
and $c_{\X}(m)$ is an element of $\R^*/\Q^*$ related to the covolume of
the image of $\widetilde{r}_{\D}$. 

Beilinson proved  this when  $\X$ is  a product of  modular curves and
$m=1$.  From  that it follows for  the product  of two elliptic curves
over $\Q$.

The purpose of this note is to directly  prove this for the product of
two elliptic curves over $\Q$. In this  situation, if $E$ and $E'$ are
isogenous,  the conjecture asserts that  the motivic cohomology is $0$
dimensional  and  hence the {\em   value}  of  $L(\X,1)$  should be  a
non-zero rational multiple of the period $c_{\X}(1)$.  Further, if $E$
and $E'$  are not isogenous,  the conjecture  asserts that the motivic
cohomology is $1$  dimensional, and hence  the value of the {\em first
  derivative}  $L'(\X,1)$  should be the   period $c_{\X}(1)$ up  to a
non-zero   rational  number.  Here this  period    is  essentially the
regulator of some element of the  motivic cohomology.  In this case we
show a      little  more  by   finding     an   explicit  element   of
$H^3_{\M}(\X,\Q(2))_{\ZZ}$ whose  regulator is equal  to  the value of
the first derivative of the $L$-function at $s=1$.  This element comes
from  a  modular parameterization  of the  elliptic   curves and is an
analogue of a cyclotomic unit.

We finally end up with an analogue of the class number formula
for real quadratic fields:

\begin{thm} Let $E$,$E'$ be non-isogenous elliptic curves over $\Q$
  corresponding modular forms  $f$ and  $g$ of  level $N_1$ and  $N_2$
  respectively. Let
$${\bf log}_q(t)=\log|qt|+\sum_{n=1}^{\infty} \log|1-q^nt|$$
where $q=e^{2\pi i z}$. Let $\xi$ be a primitive $N^{th}$ root of
unity where $N=l.c.m.(N_1,N_2)$. Then 
$$L'(H^2(E \times  E'),1)=\frac{-H(0)}{8} \mathop{\sum_{k\;mod\;N}}_{(k,N)=1}
\frac{1}{2 \pi i} \int_{X_0(N)} {\bf log}_q(\xi^k) f(q)\overline{g(q)}
\frac{dq}{q} \frac{d\bar{q}}{\bar{q}}$$
where $H(0)$ is a certain rational number corresponding to the terms
in the $L$-function for primes dividing $N$. 
\end{thm}
Stark made conjectures relating the exact values of $L$-functions of Number
fields  to regulators of units in some auxiliary number fields, and
from that point of view, this can be regarded as a special case of a
generalization of those conjectures. 

The proof follows by looking at  Ogg's \cite{og} original proof of the
Tate conjecture for products of two elliptic curves more carefully and
using Kronecker's first limit formula.

In the second part we show that the  conjecture for the product of two
elliptic   curves  implies   the conjecture   for  $H^{2n-1}_{\M}(\prod_i^n
E_i,\Q(n))_{\ZZ}$, the $n$-fold product of  elliptic curves. It appears
that one does not get any  elements excepting those induced from lower
products.

\subsection{Analogies with quadratic  extensions of $\Q$}

There is a suggestive analogy of this situation with that of quadratic
extensions of $\Q$,  which is a  special case  of $m=0$. Consider  the
group
$$\Sigma_m:=\{Ker:H^{2m+2}_{\M}(\X,\Q(m+1))_{\ZZ} 
\longrightarrow H^{2m+2}_{\M}(\X,\Q(m+1))\}$$
namely   the  group of codimensional   m   cycles supported on special
fibres.   Conjecturally   $\Sigma_m$  is finite.    When    $m=0$  and
$\X=Spec(K)$,  where  $K$ is a number  field,  $\Sigma_0$ is the class
group, which is well known to be finite.

The class number formula gives an expression for the class number
$h_K$ in terms of the Dirichlet $L$-functions associated to $K$. Let
$K^+$ denote the maximal real subfield and $h^+_K$ be the class number
of $K^+$. Let $h^*_K=h_k/h^+_K$ (so $h^+$ is the class number if $K$
is real and $h^*$ is the class number if $K$ is imaginary quadratic.) 
Then the class number formula give significantly different expressions
for $h^*$ and $h^+$. 

\subsubsection {Imaginary Quadratic Fields}

If $K$ is imaginary (so $h_K=h^*_K$ ) for  each rational prime $p$ there
is a certain  element ${\frak g}(p)$ of $K^*$  called the `Gauss  sum',
coming from a cyclotomic field  containing $K$, which has the property
that its  ideal factorization involves only the  primes lying over $p$
and does not depened on $p$. Namely
$${\frak g}(p)=\prod_{{\frak P}|p} {\frak P}^{\theta}$$
where $\theta$ is a certain element of the group ring of the integral
Galois group, the Stickleberger element, which does not depend on $p$.
The index of the ideal generated by the Stickleberger element is the
class number. Hence the special element, the Gauss sum, gives rise to
annihilators of the class group and is related to the value of the
$L$-function, though in a roundabout manner. 

Mildenhall \cite{mi} studied the group $\Sigma_1$ when $\X$ is the
self product of an elliptic curve over $\Q$. He showed that it
is {\em torsion} by constructing annihilators coming from certain
special elements of $H^3_{\M}({\bf Y},\Q(2))$ where ${\bf Y}$ is the
self product of a modular parametrisation of the elliptic curve. These
elements are analogues of Gauss sums as they too ramify at precisely
one place. However, the relation with the $L$-function is not
clear. Flach \cite{fl} studied a Selmer group associated to the
symmetric square of an elliptic curve which is conjecturally the same
as $\Sigma_1$ and did find some relation between the $L$-value and the
order of this group, though it is still not known whether this group
is finite.  

\subsubsection{Real Quadratic Fields}

Similarly, if $K$ is  real (so  $h^+=h$) and  $\chi$ is  its quadratic
character one has
$$L'(0,\chi)=h^+ \log|\epsilon| $$
where $\epsilon$ is the fundamental unit. On the other hand, one also
has the formula 
\begin{equation}
L'(0,\chi)=\log   \mathop{\prod_{k\;mod\;N}}_{(k,N)=1}
|1-\xi^k|^{-\frac{1}{2} \chi(k)} = \mathop{\sum_{k\;mod\;N}}_{(k,N)=1}
\frac{-\chi(k)}{2} \log|1-\xi^k| 
\label{cn1}
\end{equation}
where $N$  is the conductor  of $\chi$ and $\xi=e^{\frac{2\pi i}{N}}$. 
This shows that the exact value of  $L'(0,\chi)$ is the regulator of a
{\em   naturally  constructed} unit coming   form   a cyclotomic field
containg  $K$. Further, the  index of the  subgroup of the units group
generated by the cyclotomic units is the  class number. This fact is a
lot  harder to prove directly  without using the analytic class number
formula and was only done about ten years ago by Thaine \cite{th}.

Our result  can be viewed as an  analogue of the second statement \ref
{cn1} as we compute  the exact value of  the $L$-function of $E \times
E'$ in terms of the regulator of a special element coming from modular
parametrisations.   However,  in this case  a  lot less is known about
$\Sigma_1$, it is not known even whether it  is torsion and at present
it is not  clear whether one  can apply  Thaine's method to  construct
annihilators. As  far  as we are  aware, there  is  no construction of
elements   of $H^3_{\M}(E \times  E',\Q(2))_{\ZZ}$   without using the
modularity except over  a local field  by Spiess.  Such a construction
could suggest how to find an analogue of the `fundamental unit'.

Curiously it appears that  for  the imaginary quadratic and  isogenous
cases, it is easy to use the  special element to construct annihilators
but   hard to relate   to $L$-values,  while  in the  real quadratic and
non-isogenous cases, it is easy to relate the special elements to
$L$-values, but hard to construct annihilators. 

Sinnott \cite{si} has a uniform theory of  cyclotomic units and
Stickleberger elements and it is concievable that something along
those lines would generalise.

{\em Acknowledgments} The second author would like to thank S.  Bloch,
E. Ghate, K. Kimura, O. Patashnick,  K. Rogale for their comments, and
would also like to  thank  the Duke University Mathematics  department
for its hospitality when this work was started. The first author would
like to thank Queens university  for its hospitality and the Duke-IMRN
conference for providing an opportunity for this work to be completed.

\section{The Rankin-Selberg Method}

\subsection{Preliminaries}

Let $\E$ and   $\E'$  be the  two elliptic   curves  over  $\Q$.   Let
$\omega_E$ and $\omega_{E'}$  be Neron  differentials corresponding to
the  global minimal Weierstrass models.   These are defined up to $\pm
1$.   Let $f$ and  $g$ be the modular forms  of weight 2 of levels $N_1$
and $N_2$ corresponding     to   $\E$ and $\E'$   respectively.     Let
$N=l.c.m(N_1,N_2)$. We will think of $f$  and $g$ as modular forms for
$\Gamma_0(N)$. Let $\phi$ and  $\phi'$ be the modular parametrisations
from  $X_0(N)$  to $E$  and $E'$ respectively.    Define $c(\phi)$ and
$c(\phi')$ in $\Q^*$ by
$$\phi^*(\omega_E)=c(\phi)   2\pi   i   f(z) dz \;\;   \text{and}  \;\;
{\phi'}^*(\omega_{E'})=c(\phi') 2\pi i g(z) dz$$
where by $i$ we denote a choice of a $\sqrt{-1}$ that we make once and
for all. It turns out that $c(\phi)$ and $c(\phi')$ are actually in
$\ZZ\backslash \{0\}$.

Let
$$ f(z)=\sum_{n=1}^{\infty} a_n q^n  \text{    and   }
 g(z)=\sum_{n=1}^{\infty} b_n q^n $$
 be  the Fourier expansions    at $\infty$   of  $f$ and   $g$,  where
 $q=e^{2\pi  i z}$.  These modular forms   are {\em normalised} in the
 sense that they  are eigenfunctions  for all the Hecke operators for
 $p \not| N$ and for the Fricke involution, and $a_1=b_1=1$.

Let $X_0(N)$ denote  the compactification  of
the fundamental   domain for $\Gamma_0(N)$. Let

$$\delta(f,g)=f(z)\overline{g(z)}dxdy=\frac{i}{2}f(z)\overline{g(z)}dzd{\bar
  z}$$
where $z=x+iy$. Define the {\bf Petersson Inner product} by 
$$ (f,g)=\frac{1}{[\Gamma:\Gamma_0(N)]} \int_{X_0(N)} \delta(f,g)$$

We will use the following two theorems of Ogg \cite{og}.

\begin{thm}[Ogg] 

If $f$  and $g$ are normalised of  levels $N_1$ and $N_2$ respectively
and $(f,g) \neq 0$ then $f$=$g$ ( and $N_1=N_2$ ).

\end{thm}

\begin{thm}[Ogg]

If $f=\sum_{n=1}^{\infty} a_n q^n$ is a normalised cusp form of
square-free level $N$ and $p|N$ then 
$$a_p=\pm 1$$

\end{thm} 

Let $L(H^2(\E \times  \E',s))$ be the  $L$-function of the product
of the two elliptic curves. Then one has
\begin{equation}
L(H^2(E \times E'),s)=\zeta(s-1)^2 H(s-1)L_{f,g}(s-1)
\label{lh2}
\end{equation}
where   $\zeta(s)$ is the Riemann Zeta function, $\zeta_N(s)$ is the
same function with the primes dividing $N$ removed, $H(s)$ is a
polynomial in $p^{-s}$ coming from the primes dividing $N$ and
$$L_{f,g}(s)=\zeta_N(2s)   \sum_{n=1}^{\infty}  \frac{a_n
{\overline b_n}}{n^{s+1}}$$

\subsection{Rankin-Selberg Convolution}

We use the Rankin-Selberg convolution to get an integral
representation of the $L$-function. 

Let $f$,$g$ and $z$ be as above. Then

$$\int_{\frac{-1}{2}}^{\frac{1}{2}}      f(z)\overline    {g(z)}   dx=
\sum_{n=1}^{\infty} a_n \overline {b_n} e^{-4\pi n y} $$
and so, integrating over the region $y > 0$ as well ,

$$  \int_{|x|\leq   \frac{1}{2}} f(z) \overline{g(z)}y^{s-1}  dx   dy =
(4\pi)^{-s}\Gamma(s)\sum_{n=1}^{\infty}\frac{a_n
\overline{b_n}}{n^{-s}}$$
Replacing $s$ by $s+1$ we have 
\begin{equation}
(4\pi)^{-s-1} \Gamma(s+1) \sum_{n=1}^{\infty}  \frac{a_n  {\overline
b_{n}}}{n^{-(s+1)}}=\int_{|x|\leq \frac{1}{2}} y^s \delta(f,g)
\label{rs1}
\end{equation}

The  region $|x|\leq \frac{1}{2}$  is  the fundamental domain for  the
stabilizer of the cusp $\infty$

$$\Gamma_{\infty}=\Gamma_0(N)_{\infty}=\{  \pm  \begin{pmatrix}  1&n\\0&1
\end{pmatrix} |n \in \ZZ \}$$
so one has 
$$\int_{|x|\leq      \frac{1}{2}}   y^s       \delta(f,g)=\int_{{\frak
h}/\Gamma_{\infty}} y^s  \delta(f,g)=\int_{X_0(N)} E_{\infty}^{N}(z,s)
\delta(f,g)$$
where $E_{\infty}^{N}(z,s)$ is the {\bf Eisenstein Series} 
$$ E_{\infty}^{N}(z,s)=\sum_{\gamma   \in  \Gamma_0(N)/\Gamma_{\infty}}
(Im(\gamma z))^s = 1+ \mathop{\sum_{m >0}}_{(mN,n)=1} \frac{y^s}{|mNz+n|^{2s}}$$
Let 
$$\zeta_N(s)=\mathop{\sum_{n>0}}_{(n,N)=1} n^{-s}=\prod_{p \not \; | N}(1-p^{-s})^{-1}$$
One then has 
$$2\zeta_N(2s) E^{N}_{\infty}(z,s)=\sum_{d | N} \frac{\mu(d)}{d^s} E_{\infty}(\frac{Nz}{d},s)$$
where $\mu(d)$ is the M\"{o}bius function and 
$$E_{\infty}(z,s)=\sum_{m,n}{}' \frac{y^{s}}{|mz+n|^{2s}}$$
Let  
$$L_{f,g}(s)=\zeta_N(2s)\sum_{n=1}^{\infty} a_n\overline{b_n}n^{-(s+1)} $$
Subsituting this into the  equation \ref{rs1} we get
\begin{equation}
2(4\pi)^{-s-1}\Gamma(s+1)L_{f,g}(s)=\sum_{d|N}\frac{\mu(d)}{d^{2s}}
\int_{X_0(N)} \delta(f,g) E_{\infty}(\frac{Nz}{d},s)
\label{rs2}
\end{equation}

\subsubsection{ The Epstein-Zeta Function}

The function 
$$E_{\infty}(z,s)=  \sum_{\gamma \in  \Gamma/\Gamma_{\infty}}
 Im(\gamma   z)^s= \sum_{m,n} {}'
\frac{y^s}{|mz+n|^{2s}}$$
appears  in many different guises  and is sometimes  known as the {\bf
Epstein Zeta Function} or an {\bf Eisenstein-Kronecker-Lerch Series}.
It converges for $Re(s)>1$  and has  a meromorphic continuation to  the
entire complex  plane. Further, the function
$$E^*_{\infty}(z,s)=(\frac{1}{\pi})^{s} \Gamma(s)E_{\infty}(z,s)$$
satisfies the functional equation 
$$E^*_{\infty}(z,s)=E^*_{\infty}(z,1-s)$$
and has a simple pole  with residue $1$  at $s=1$ {\em independent  of
  $z$}. A good reference  for all these facts  is Lang's {\em Elliptic
  Functions} \cite{la}.

\subsubsection{An Integral Representation}

Multiplying by  $(\frac{N}{\pi})^2 \Gamma(s)$  and using the  function
$E^*_{\infty}(z,s)$ in equation \ref{rs2} above gives
$$\Phi(s):=(\frac{2\pi}{\sqrt{N}})^{-2s}\Gamma(s)\Gamma(s+1)L_{f,g}(s)        =
2\pi\sum_{d|N}\frac{\mu(d)}{d^{s}}       \int_{X_0(N)}     \delta(f,g)
E^*_{\infty}(\frac{Nz}{d},s)  $$
Since the residue of $E_{\infty}^*(z,s))$ at $s=1$ is $1$ {\em independent
of $z$} and it is a simple pole one sees that the residue of $\Phi(s)$
at $s=1$ is a constant times $(f,g)$. From that one has

\begin{thm}[Rankin]

$$L_{f,g}(s)=\zeta_{N}(2s)\sum_{n=1}^{\infty}a_n\overline{b_n}n^{-(s+1)}$$
is entire if $(f,g)=0$ and is entire except for a simple pole at $s=1$ if 
$(f,g)\neq 0$. In this case the residue is a rational number times $(f,g)$.

\end{thm}

\section{The Tate Conjecture for $E \times E'$}

\noindent The  Tate conjecture amounts to  the following two statements: 

\begin{itemize}

\item $L(H^2(E \times E'),s)$ has a pole of order $3$ at $s=2$
  when $\E$ and $\E'$ are isogenous. 
  
\item $L(H^2(E \times  E'),s)$ has a pole of  order $2$ at  $s=2$ when
  $\E$ and $\E'$ are {\bf not} isogenous

\end{itemize}

Note that  it does not matter if  $E$ has complex multiplication as we
are looking at the  cycles defined  over  $\Q$ and the extra  cycle is
only defined over the field of CM.

Since  $\zeta(s)$     has     a   simple      pole   at   $s=1$, using
\ref{lh2} the   conjecture reduces to the statements that

\begin{itemize}

\item  $L_{f,g}(s)$ has  a simple pole at $s=1$ if $E$ and $E'$ are isogenous.

\item   $L_{f,g}(s)$  is    holomorphic    and   non-vanishing at      $s=1$  
if  $E$ and $E'$ are {\bf not} isogenous.

\end{itemize}

From  Rankin's theorem one has that  $L_{f,g}(s)$ has a simple pole at
$s=1$ when $E$ and $E'$ are isogenous as $(f,g)  \neq 0$.  When $E$ and
$E'$ are not isogenous, $(f,g)=0$ hence there is no pole at $s=1$.

To complete the proof of the Tate  conjecture, we use the following
theorem of Ogg's \cite{og}

\begin{thm}[Ogg]
$L_{f,g}(1)   \neq    0$ if  $(f,g)=0$.  
\label{og}
\end{thm}

The  proof   of  this theorem  is   by  using  the  Euler product  for
$L_{f,g}$ to construct a Dirichlet   series with positive   real
coefficients which does  not have a pole  contradicting the  fact that
such a  Dirichlet series  has a pole on   the real point of  the
critical line. Details can be found in Ogg's paper \cite{og}.

\section{Beilinson's Conjecture for $E \times E'$}

To  verify  Beilinson's  conjecture, we    first have  to  get    some
understanding of what the integral  motivic cohomology groups are.  In
the following   sections, we describe the   group $H^3_{\M}(X,\Q(2))$,
where $X$ is  a  surface  defined over   $\Q$.  The integral   motivic
cohomology  group $H^3_{\M}(X,\Q(2))_{\ZZ}$  is  a certain subgroup of
this  group, first described   conditionally by Beilinson, though more
recently unconditionally  by Scholl \cite{sc}. 

\subsection { Elements of $H^3_{\M}(X,\Q(2))$ }

Let $X$  be a surface  defined  over $\Q$.  The  group  $H^3_{\mathcal
  M}(X,\Q(2))$ has several  different descriptions: First, in terms of
a   graded  piece of   $K_1(X) $, second,   as  the  higher Chow group
$CH^2(X,1)$ and finally, as  the $K$-cohomology group $H^1(X,{\mathcal
  K}_2)$.   From   the third    description  and the   Gersten-Quillen
resolution, an element of the group is represented by a formal sum
$$\sum  ({\mathcal C},f)$$
where $\C$ are curves on $X$  and $f$ are functions on these
curves subject to the cocycle condition
$$\sum div(f)=0$$
This  is a generalization of   the fact that   elements  of $F^*$  are
elements of $K_1$ of a number field $F$.  

\subsubsection{ Construction of the elements on products of curves }

We construct  some elements of  the group $H^3_{\mathcal M}(X,\Q(2))$ when
$X$ is the  self product of a  curve  $\C$. We  use  a construction of
Bloch's which was generalized by Beilinson \cite{be}.

Let $\C$ be a  curve which contains a  set $\CC$ such that any divisor
of degree $0$ supported  on $\CC$ is torsion in  the jacobian of $\C$. 
Let $Y=\C \backslash \CC$.  The condition above  can be stated as  the
statement that the exact sequence
\begin{equation} 
0 \rightarrow H^1_{\mathcal M}(\C,\Q(1)) \rightarrow   H^1_{\mathcal
  M}(Y,\Q(1)) \stackrel{\partial}{\rightarrow}  (H^0_{\M}(\CC,\Q(0)))_0
\rightarrow 0
\label{ex1}
\end{equation}
{\em splits} as a sequence  of motivic cohomology  groups as the class
in the Chow group of a divisor of degree $0$ can be described in terms
of   such an  extension. The  groups  $H^1_{\mathcal M}(\C,\Q(1))$ and
$H^1_{\mathcal M}(Y,\Q(1))$   are  simply $\OO^*_{\C}  \otimes \Q$  and
$\OO^*_{Y} \otimes \Q$ respectively and
$$(H^0_{\M}(\CC,Q(0)))_0=((z_1,z_2,...,) \in  \bigoplus_{s
\in \CC} \ZZ) \text { such that } \sum z_i=0$$
Splitting   is the statement  that for   a divisor  $D$ of degree  $0$
supported   on $\CC$, there is    a  canonical choice  of a   function
$\epsilon(D)$ whose divisor is $D$.

One then has the following Lemma:

\begin{lem}
The sequence 
\begin{equation}
0 \rightarrow       H^3_{\mathcal      M}({\C}^2,\Q(2))  \rightarrow
H^3_{\mathcal       M}({\C}^2      \backslash     {\CC}^2      ,\Q(2))
\stackrel{\partial}{\rightarrow}           (H^0_{\M}({\CC}^2,\Q(0)))_0
\rightarrow 0
\end{equation}
splits as a sequence of motivic cohomology groups. 
\end{lem}

\begin{proof}
  
  The idea is   to use  the   splitting in the first  case,   equation
  \ref{ex1}, to split it in this case.  We need to produce a canonical
  element  of    $H^3_{\M}(\C^2  \backslash \CC^2   ,\Q(2))$  for  any
  $0$-cycle of  degree $0$ supported  on  $\CC^2$.  It suffices to  do
  this for cycles of the form $D=(P_1,P_2)-(Q_1,Q_2)$.  We use a trick
  familiar to all first year students of calculus.
$$(P_1,P_2)-(Q_1,Q_2)=(P_1,P_2)-(Q_1,P_2)+(Q_1,P_2)-(Q_1,Q_2)$$
For the   pair $(P_1,P_2)-(Q_1,P_2)$ we  take  the element $(Y \times
P_2,   \epsilon  (P_1-Q_1))$     and   similarly,    for   the    pair
$(Q_1,P_2)-(Q_1,Q_2)$ we  take the  element $(Q_1 \times  Y, \epsilon
(P_2-Q_2))$. The sum of these two elements  give the canonical lift of
$D$.  A  different way of splitting $D$  gives the same element as two
such liftings  differ by something coming from  the tame symbol, which
is a coboundary. Since any  divisor of degree $0$ can  be written as a
sum of such $D$, this gives a splitting.

\end{proof}

For a divisor $D$ of  degree $0$ on  $\CC^2$, let $\epsilon^2(D)$ denote
the  lifting. Now suppose one has a  map $\Psi:X  \rightarrow X^2$ such
that
\begin{itemize}
\item $\Psi(Y) \subset X^2 \backslash \CC^2$
\item $\Psi(\CC)  \subset \CC^2$
\end{itemize}
\noindent One then has  an induced pushforward, 
$$\Psi_*:H^1_{{\M}}(Y,\Q(1))   \longrightarrow   H^3_{\M}(X^2
\backslash \CC^2,\Q(2))$$
so if $D$ is a divisor of degree $0$ on $\CC$, one gets an element 
$$\Psi_*(\epsilon(D)) \in H^3_{\M}(X^2 \backslash \CC^2,\Q(2))$$
The simplest example of such a map $\Psi$ is the diagonal embedding.

Now if  $D$ is a divisor  of degree $0$  on $\CC$ one gets two elements,
$\epsilon^2(\Psi(D))$  and  $\Psi_*(\epsilon(D))$   in     $H^3_{\M}(X^2
\backslash            \CC^2,\Q(2))$.         The                 element
$\chi(D)=\epsilon^2(\Psi(D))-\Psi_*(\epsilon(D))$ satisfies the condition
$$\partial(\epsilon^2(\Psi(D))-\Psi_*(\epsilon(D)))=\Psi(D)-\Psi(D)=0$$
Hence it lifts to give an element of $H^3_{\M}(X^2,\Q(2))$ !  

\subsubsection{ Products of modular  elliptic curves}

Returning to our case,  let  $N=l.c.m(N_1,N_2)$ where $N_1$ and  $N_2$
are the  conductors of $E$ and  $E'$ respectively, and  let $X_0(N)$ be
the modular curve of level $N$. One then has modular parametrisations,
$$\phi:X_0(N) \longrightarrow E \text { and } \phi':X_0(N) 
\longrightarrow E'$$  
We  take $\C=X_0(N)$. By the Manin-Drinfel'd   theorem, any divisor of
degree  $0$ supported on the set  of cusps is  torsion, so we can take
this set as the set $\CC$. Applying the  above lemma, we can construct
elements of    $H^3_{\M}(X_0(N)^2,\Q(2))$   and   using  the   modular
parametrizations,  we   can push these   elements  down to $H^3_{\M}(E
\times E',\Q(2))$.

\subsubsection{ Remarks on Integrality}

Beilinson's conjecture is about the {\em integral} motivic cohomology.
This  is a subspace  of  the motivic  cohomology which  was originally
defined to be the image of the motivic cohomology  of a regular proper
model, if it exists. Scholl \cite{sc} gave an unconditional definition
of  this  subspace,  denoted by $H^3_{\M}(X,\Q(2))_{\ZZ}$,  using  De
Jong's theory of alterations. This is the analogue of $\OO_K^*$ as
opposed to $K^*$ where $K$ is a number field.

In general  the elements we construct using  the above method  do {\em
  not}   lie in  the  integral  motivic  cohomology.  However,  Scholl
\cite{sc}  showed   that   the   projection    of  the   elements   of
$H^3_{\M}(X_0(N) \times X_0(N),\Q(2))$  onto  the  $H^3_{\M}(E  \times
E',\Q(2))$ lie in $H^3_{\M}(E \times  E',\Q(2))_{\ZZ}$ if $E$ and $E'$
are  not  isogenous.  Further Harris  and Scholl  \cite{sc2} show that
this  subspace  of    $H^3_{\M}(E  \times E',\Q(2))_{\ZZ}$   is   zero
dimensional if $E$ and $E'$ are  isogenous and at most one dimensional
if $E$ and $E'$ are not isogenous. It is not known if  $H^3_{\M}(E
\times E',\Q(2))_{\ZZ}$ is even finitely generated.

To  show that it is at  least  one, we will use  the  fact that if the
regulator of an element is non-zero,  then the element cannot be zero. 
This regulator turns out to be the  value $L'_{f,g}(0)$.  From Theorem
\ref{og} and  the functional equation,  we know  that this is non-zero
and we can conclude  that part of the  Beilinson conjecture.   Towards
that end, in the next section we will describe this regulator.

If one does not require   integrality, Flach \cite{fl} and  Mildenhall
\cite{mi} have shown independently  that the group $H^3_{\M}(E  \times
E',\Q(2))$ is infinitely generated.

\subsection{ Regulator Maps and the Real Deligne Cohomology}

Let $\X$ be a smooth projective variety over $\Q$, $i$ an even integer
and  $m=i/2$.  In this section     we describe what the   Real
Deligne  cohomology   groups   are and  explain    what  the constants
$c_{\X}(m)$   are.  Details    can  be found    in   the  articles  of
Esnault-Viehweg and Schneider in \cite{ssr}.
 
The {\bf  regulator map} $\widetilde{r}_{\D}$  has two components,
$\widetilde{r}_{\D}=r_{\D}\oplus z_{\D}$, where 
\begin{itemize}
\item $z_{\D}$ is a cycle class map induced by the usual cycle class map
  to De Rham cohomology,
$$z_{\D}:B^{m}(\X)_{\Q} \longrightarrow  H^{i+1}_{\D}
(\X_{/\R},\R(m+1))$$
\item $r_{\D}$ is a higher cycle class map, generalizing Dirichlet's
regulator map for units,
$$r_{\D}:H^{i+1}_{\M}(\X,\Q(m+1))_{\ZZ}  \longrightarrow H^{i+1}_{\D}   (\X_{/\R},\R(m+1))$$
\end{itemize}
Here   $H^{i+1}_{\D}(\X_{/\R},\R(m+1))$  is   the  {\bf  Real  Deligne
  cohomology} described below.

\subsubsection {Real Deligne Cohomology}

The Real Deligne  cohomology is a  real vector space  generalising the
vector space  $\R^{r_1+r_2}$, which  appears as  the target  space for
Dirichlet's regulator map. While the precise  defininition of the Real
Deligne   cohomology is  a   little involved, there    are  two  key
properties:
\begin{itemize}

\item (i) There is an exact sequence, 
$$0\rightarrow F^{m+1}H^{i}_{DR}(\X_{/\R}) \rightarrow H^i_{B}(\X({\Bbb
  C}),\R(m))^{(-1)^{m}} \rightarrow$$
\begin{equation}
\rightarrow H^{i+1}_{\D}(\X_{/\R},\R(m+1))\rightarrow 0
\label{ex2}
\end{equation}
where the $-1$ indicates that it is a $-1$ eigenspace for the
involution induced by complex conjugation on the complex manifold
$X({\Bbb C})$, $H_B$ is the Betti ( singular ) cohomology and $H_{DR}$
is the algebraic De Rham cohomology. 

\item (ii) The dimension is related to the order of vanishing of
  $L$-functions:
$$dim_{\R} H^{i+1}_{\D}(\X_{/\R},\R(m+1))=ord_{s=m} L(H^i(\X,s))-ord_{s=m+1}L(H^i(\X,s))$$

\end{itemize}

\subsubsection{The period $c_{\X}(m)$}

From \ref{ex2} there is an isomorphism of one dimensional vector
spaces 
$$det(F^{m+1}H^{i}_{DR}(\X_{/\R})) \otimes det (H^{i+1}_{\D}(\X_{/\R},\R(m+1)))
\simeq det (H^i_{B}(\X({\Bbb C}),\R(m))^{(-1)^{m}}$$
As   $\X$ is defined   over $\Q$, $F^{m+1}H^{i}_{DR}(\X_{/\R})$ has a
rational structure  coming  from  the algebraic De  Rham   cohomology. 
$H^{i}_{B}(\X({\Bbb C}),\R(m))^{(-1)^{m}}$  has  an  obvious  rational
structure.  Part   (i)   of  the Beilinson   conjecture  asserts  that
$Im(\widetilde{r}_{\D})$ gives a   rational structure  on  the Deligne
cohomology    $H^{i+1}_{\D}(\X_{/\R},\R(m+1))$.  So all the vector spaces
involved have rational structures, at least conjecturally. $c_{\X}(m)$
is the determinant of the isomorphism above  computed with respect
to these $\Q$-structures. It is an element of $\R^*/\Q^*$.

\subsubsection{Explicit formulae for the Regulator map}

The regulator  map has the   following explicit  description as a
current on $(m,m)$ forms:

\begin{itemize}
\item On $B^m(\X)_{\Q}$ it   is given   by a  current   of
  integration. If  $\Z$  is an  element  of  $B^m$  and $\omega$  is a
  $(dim(X)-m,dim(X)-m)$ form in $H^i_{DR}(E \times E')$ then
\begin{equation}
(z_{\D}(\Z),\omega) := (\frac{1}{2 \pi i})^{m} \int_{\Z} \omega
\label{reg1}
\end{equation}

\item  On the motivic cohomology side.  If $\sum (\C,f)$ is an element
  of   $H^{i+1}_{\M}(E   \times   E',\Q(m+1))$   and  $\omega$   is  a
  $(dim(X)-m,dim(X)-m)$ form then
\begin{equation}
(r_{\D}(\sum (\C,f)),\omega) := (\frac{1}{2\pi i})^{m} \sum \int_{\C}
\log|f| \omega
\label{reg2}
\end{equation}
\end{itemize}

This regulator map is  conjecturally  injective. However, it  is clear
that  if $(r_{\D}(\sum (\C,f),\omega)  \neq 0$ for some $\omega$, then
the element is non-trivial in $H^{i+1}_{\M}(E \times E',\Q(m))$.

We are interested in the case $\X=E \times E'$ and $i=2$ so $m=1$. In
this case the corresponding Deligne cohomology is $3$
dimensional. However, it further breaks down according to the motivic
decomposition described in the next section.  

\section{Calculation of $L$-values } 

\subsection{Motivic Decomposition}

To  compute  the $L$-values   we  first observe  that  we   can use  the
K\"{u}nneth formula to get a decomposition of the motive $H^2(E \times
E')$.

The motive $H^2(E \times E')$ splits up in to 4 or 3
submotives depending on whether $E$ is isogenous or not to $E'$.
From the K\"{u}nneth formula one has 
\begin{equation}
H^2(E \times E')=H^2(E) \oplus H^2(E') \oplus H^1(E)\otimes H^1(E')
\label{ku}
\end{equation}
If $E \simeq E'$ then the motive $H^1(E)\otimes H^1(E')$ further splits
up   in       to 
\begin{equation}     
H^1(E)\otimes H^1(E)=\Lambda^2H^1(E)    \oplus     Sym^2    H^1(E) 
\label{ku2}
\end{equation}
and $\Lambda^2H^1(E)\simeq H^2(E)$.

As  this   decomposition is    at   the level  of motives,    all  the
corresponding objects     such as  $L$-functions   and    the  constants
$c_{\X}(1)$ also decompose and we can reduce  the problem to verifying
each case individually.

The  case  of real   interest  to us is   that   of $L(H^1(E)  \otimes
H^1(E'),s)$ when  $E$ and $E'$ are  {\em not}  isogenous, as the other
cases, namely $L(H^2(E))$,  $L(H^2(E'))$ and $L(Sym^2(H^1(E)))$ either
reduce to the cases of fields,  as in the  first two cases, or have been
treated in detail elsewhere \cite{fl}, as in the third. However, for
completeness we will describe them. 

\subsubsection {The motives $H^2(E)$ and $H^2(E')$}

In either  case  here the Deligne   cohomology is  $1$ dimensional  as
$H^2(E)$  is $1$  dimensional, complex   conjugation acts by  $-1$  and
$F^2H^2_{DR}=0$. From the exact sequence \ref{ex2} we have
$c_{H^2(E)}(1)$ is given by 
$$\frac{1}{2 \pi i} \int_E \alpha$$
where $\alpha$ is a rational De Rham  cohomology class. Such a form is
obtained by $\omega_E \wedge \eta_E$ where $\omega_E$ is the canonical
differential  and  $\eta_E$   is the  form  $\frac{XdX}{Y}$.  The {\em
  Legendre relation} shows
\begin{equation}
\frac{1}{2 \pi i} \int_E \omega_E \wedge \eta_E =1
\label{h2e}
\end{equation}
so it is rational.

The   $L$-function is   $L(H^2(E),1)=\zeta(s-1)$ and so
$$L(H^2(E),1)=\zeta(0)=\frac{1}{2} \sim_{\Q^*}  c_{H^2(E)}(1)$$

\subsubsection {Functional Equation} 

For   the other two  cases, we   need the  functional equation  of the
$L$-function $L_{f,g}(s)$  to compute the value  at $s=1$.  Here we make
the  further assumption   that   $N_1$     and  $N_2$  and       hence
$N=l.c.m.(N_1,N_2)$ are square-free. Then one has

\begin{thm}[Ogg] 

Let $f$ and $g$  be  normalized cusp  forms of square-free  level
$N_1,N_2$; Let $N=l.c.m.(N_1,N_2)$ and $M=g.c.d.(N_1,N_2)$. For $p|M$,
let
$$c_p=a_pb_p=\pm 1$$
Recall
$$\Phi(s):=(\frac{2\pi}{\sqrt{N}})^{-2s}\Gamma(s)\Gamma(s+1)L_{f,g}(s)$$
and  set 
$$\Phi^{+}(s)=\Phi(s)A(s)$$
where 
$$A(s)=\prod_{p|M}(1-c_pp^{-s})^{-1}$$
Then 
$$\Phi^{+}(s)=\Phi^{+}(1-s)$$

\end{thm}

A generalization   of this theorem  for arbitrary  $N$ was   proved by
\cite{li}.  

\noindent The $L$-function of $H^1(E)\otimes H^1(E')$ is
\begin{equation}
L(H^1(E) \otimes H^1(E'),s)=H(s-1)L_{f,g}(s-1)
\label{lh1}
\end{equation}
where $H(s)$ is a term depending on primes dividing $N$. In
general it can be quite complicated, but the following ad hoc
definition, due to Ogg \cite{og}, seems to make the formulae cleaner.

If $p|M$, then the factor is 
$$\frac{1}{(1-c_p p^{-s})(1-c_p p^{-(s+1)})}$$ 
while if $p|N_1$, $p \not | N_2$, it is 
$$\frac{1}{(1-a_p b_p p^{-(s+1)} + p^{-1-2s})}$$
at least when $E$ has multiplicative reduction at $p$. 

\subsection{ Assume $E$  isogenous to $E'$} 

In this case, the motivic decomposition shows that we have to
understand the case of $Sym^2(E)$. This is a {\em critical} motive in
the sense of Deligne \cite{dl} as the Deligne cohomology
vanishes. The period turns out to be \cite{dl}
\begin{equation}
c_{Sym^2(E)}(1)=\frac{1}{2 \pi i} \int_{E({\Bbb C})} \omega_{E} \wedge 
\overline{\omega_{E}}
\label{csym}
\end{equation}
where $\omega_E$ is the canonical differential.

Here $N_1=N_2=N$ and $f=g$ so $c_p=a_p^2=1$  for all $p|N$. Let $m$ be
the  number    of   primes dividing     $N$.   From  Rankin's theorem,
$L_{f,g}(s)=L_{f,f}(s)$  has a   pole  at  $s=1$.  From  the   motivic
decomposition, we have
$$L(H^1(E)\otimes H^1(E),s)=L(H^2(E),s)L(Sym^2(E),s)$$
so from \ref{lh1} above, we have 
\begin{equation}
L(Sym^2(E),s)=\frac{H(s-1)L_{f,f}(s-1)}{\zeta(s-1)}
\label{ls2} 
\end{equation}
From the functional equation we get 
$$\Phi^{+}(0)=\Gamma(0) \Gamma(1) L_{f,f}(0) A(0) = A(0) \Phi(0)$$
and $A$ has  a pole of order $m$ at $0$.

From the integral expression for $\Phi(s)$ we have, 
$$\Phi(0)=\lim_{s\rightarrow 0}  2 \pi \sum_{d|N}   \frac{\mu(d)}{d^s}
\int_{X_0(N)} \delta(f,f) E^*_{\infty}(\frac{Nz}{d},s)$$
Since  $E^*_{\infty}(\frac{Nz}{d},0)$  has a simple pole  with residue
$1$ at $s=0$ and $$\sum_{d|N} \mu(d)=\prod_{p|N}(1-1)=\frac{1}{A(0)}$$
one has
\begin{equation}
L_{f,f}(0)=\frac{\Phi(0)}{\Gamma(0)}=2 \pi \lim_{s      \rightarrow     0}
\frac{1}{s} \int_{X_0(N)} \delta(f,f) E^*_{\infty}(\frac{Nz}{d},s)  = 2
\pi (f,f)[\Gamma:\Gamma_0(N)]
\label{lff1}
\end{equation}
From the relation between the canonical differential and the modular
form one has 
$$(f,f)=\frac{i deg(\phi)}{8 \pi^2 c(\phi)^2} \int_{E({\Bbb C})}
\omega_E \wedge \overline{\omega_{E'}}$$ 
Using this in the formula \ref{lff1} and the expression \ref{csym} one gets 
\begin{equation}
L_{f,f}(0)=\frac{-deg(\phi)[\Gamma:\Gamma_0(N)]}{2c(\phi)^2} \frac{1}{2 \pi i} \int_{E({\Bbb C})}
\omega_E \wedge \overline{\omega_{E'}} \sim_{\Q^*} c_{Sym^2(E)}(1)
\label{lff2}
\end{equation}
Finally we have, from   \ref{lh1} and $\zeta(0)=1/2$,
\begin{equation}
L(Sym^2{E},1)=\frac{-deg(\phi)H(0)[\Gamma:\Gamma_0(N)]}{c(\phi)^2}
\frac{1}{2     \pi    i}  \int_{E({\Bbb C})}     \omega_E  \wedge
\overline{\omega_{E'}} \sim_{\Q^*} c_{Sym^2(E)}(1)
\label{lsym1}
\end{equation}
which is precisely what the conjecture predicts.

\subsection{ Assume $E$  is not  isogenous to $E'$} 

\vspace{\baselineskip}

In this case the $L$-function is 
\begin{equation}
L(H^1(E)\otimes H^1(E'),s)=H(s-1)L_{f,g}(s-1)
\label{lfg}
\end{equation}
The Deligne cohomology is one  dimensional and it turns out  \cite{be}
that   the  period $c_{H^1  \otimes H^1}(1)$  is  the regulator  of an
element of the motivic cohomology.  So we have to  show that the value
of  $L'(H^1(E)\otimes H^1(E'),1)$ is rational  up  to the regulator of
such an element. From \ref{og}, we know this value is non-zero, so
that will also show that the element is non-zero. 

We  assume for simplicity that $N_1$ and  $N_2$ are coprime so
$\Phi^{+}(s)=\Phi(s)$ and hence $\Phi(s)=\Phi(1-s)$. Since 
$$\Phi(0)=\Gamma(0)\Gamma(1)L_{f,g}(0)$$
and  $\Gamma(0)$ has a simple pole and $L_{f,g}(0)=0$, we have
\begin{equation}
L'_{f,g}(0)=\Phi(0)
\label{lfgp}
\end{equation}

\subsubsection { The Kronecker Limit Formula}

To compute $\Phi(0)$ we     need  {\bf Kronecker's  First    Limit
  Formula}.    This  allows us  to compute  the  constant term  in the
Laurent series expansion of $E_{\infty}(z,s)$.

\begin{thm}[Kronecker] Let 
$$E_{\infty}(z,s)=\sum_{m,n}{}' \frac{y^s}{|mz+n|^{2s}}$$
and let 
$$\eta(z)=q^{1/24}\prod_{n=1}^{\infty}(1-q^n)$$
where $q=e^{2\pi i z}$. Let $\gamma$ be Euler's constant. Then 
\begin{equation}
E_{\infty}(z,s)=  \frac{\pi}{s-1}- \pi \log y +
2\pi (\gamma-\log 2) -4\pi \log |\eta(z)| + O(s-1)
\label{klf}
\end{equation}

\end{thm}
\begin{proof}
The proof of this theorem can be found in Lang \cite{la}
\end{proof}

\noindent Recall 
$$E^*_{\infty}(z,s)=\frac{1}{\pi^s}\Gamma(s) E_{\infty}(z,s)$$
From the functional equation, we have 
$$E^*_{\infty}(\frac{Nz}{d},s)=E^*_{\infty}(\frac{Nz}{d},1-s)$$
Combining this with the limit formula \ref{klf}, we have 
$$\lim_{s \rightarrow 0} E_{\infty}^*(\frac{Nz}{d},s)=\lim_{s
  \rightarrow 0} E_{\infty}^*(\frac{Nz}{d},1-s)$$
$$=\frac{1}{-s}-  \log \frac {Ny}{d} + 2(\gamma-\log 2)
-4\log |\eta(\frac{Nz}{d})| + O(-s)$$
From this it follows that 
$$\lim_{s \rightarrow 0} \sum_{d|N} \frac{\mu(d)}{d^s}
E_{\infty}^*(\frac{Nz}{d},s)=\lim_{s \rightarrow 0} \sum_{d|N}
 \frac{\mu(d)}{d^s} E_{\infty}^*(\frac{Nz}{d},1-s)$$
$$=\sum_{d|N} -4\log |\eta(\frac{Nz}{d})|$$
As all the other terms vanish from the facts that, for $N>1$, 
$$\sum_{d|N} \mu(d)=0 \text{  and } \prod_{d|N} (\frac{N}{d})^{\mu(d)}=1$$

\noindent Recall, 
$$\Phi(s)=2 \pi \sum_{d|N} \frac{\mu(d)}{d^s} \int_{X_0(N)}
E^*_{\infty}(\frac{Nz}{d},s) \delta(f,g)$$
Using the above result we have 
\begin{equation}
 \Phi(0)=2\pi \int_{X_0(N)}  \sum_{d|N} -4 \log |\eta(\frac{Nz}{d})|
\delta(f,g)
\label{phi0}
\end{equation}

Let $\Delta(z)$ be the usual cusp form of weight 12 for
$SL_2(\ZZ)$. We have $\Delta(z)=\eta(z)^{24}$. Define   
$$\Delta_N(z):=\prod_{d|N} \Delta(\frac{Nz}{d})^{\mu(d)}$$
This is a modular {\em unit}  as $\sum_{d|N} \mu(d)=0$ and its divisor
is supported on the cusps.

Thus, the Eisenstein series tends to $\frac{1}{24}\log|\Delta_N(z)|$ and one has
\begin{equation}
\Phi(0)=\frac{-\pi}{3} \int_{X_0(N)} \log|\Delta_{N}(z)| \delta(f,g) 
\label{eis1}
\end{equation}

\subsubsection{ A natural element of $H^3_{\M}(X_0(N) \times X_0(N),\Q(2))$}

Let $D$ be the  divisor of $\Delta_N$.  Then $\epsilon(D)=\Delta_N$ and
one can consider   the  element $\chi(D)$ in  $H^3_{\M}(X_0(N)  \times
X_0(N),\Q(2))$ constructed in section 5.1. There is a map
$$\psi: X_0(N) \times X_0(N) \rightarrow E \times E'$$ 
induced by  the modular parametrizations $\phi$  and $\phi'$  and one
has the element $(\psi_*(\chi(D))$ in $H^3_{\M}(E \times E',\Q(2))$.

Let $\omega_E$ and $\omega_{E'}$ be the canonical differentials on $E$
and  $E'$  respectively.  Then  one  has   $\phi^*(\omega_E)=2 \pi   i
c(\phi)f(z)dz$ and ${\phi'}^*(\omega_{E'})=2 \pi i c(\phi')g(z)dz$.
So from that, 
$$\delta(f,g)=\frac{deg(\psi)}{8 \pi^2 i c(\phi) c(\phi')} \omega_E \wedge
\overline{\omega_{E'}}$$
where $deg(\psi)$ denotes the degree of the map $\psi|_{diagonal}$.

Therefore one has 

\begin{equation}
L'_{f,g}(0)=\frac{\pi}{3} \int_{X_0(N)} \log|\Delta_{N}(z)|
\delta(f,g)
\label{lv1}
\end{equation} 

\begin{equation}
=\frac{2\pi^2 i}{3} \frac{1}{2\pi i} \int_{X_0(N)} \log|\Delta_{N}(z)|
\delta(f,g)
\label{lv2}
\end{equation}

\begin{equation}
=\frac{2\pi^2 i}{3} \frac{1}{2 \pi i} \int_{\psi_{*}(X_0(N))}
\log|\psi_*(\Delta_{N}(z))| \frac{deg(\psi)}{8 \pi^2 i c(\phi)
  c(\phi')} \omega_E \wedge \overline{\omega_{\E'}}
\label{lv3}
\end{equation}
Using the explicit description of the regulator map \ref{reg2} we have
\begin{equation}
L'_{f,g}(0)=\frac{-deg(\psi)}{12 c(\phi) c(\phi')}(r_{\D}({\psi}_*(\chi(D))),\omega_E           \wedge            \overline
{\omega_{E'}})
\label{lv4}
\end{equation}
Since $L'(H^1(E) \otimes H^1(E'),1))=H(0)L'_{f,g}(0)$ we have  
$$L'(H^1(E)  \otimes H^1(E'),1))=\frac{-deg(\psi) H(0)}{12 c(\phi) c(\phi')} (r_{\D}( {\psi}_*(\chi(D))   ),
\omega_E \wedge \overline {\omega_{E'}})$$
$$\sim_{\Q^*} c_{H^1(E)\otimes H^1(E')}(1)$$

\subsection{A `class number formula'}

To add credence to our claim that the element we have is an analogue
of a cyclotomic unit, we show there is a `class number formula'
analogous to the expression \ref{cn1} for $\zeta'_K(0)$. 

We have   the following  curious   product formula for    our function
$\Delta_N(z)$ which can be found in Asai \cite{as}.
$$ \Delta_N(z)=\sum_{d|N}
\Delta(\frac{Nz}{d})^{\mu(d)}=q^{\phi(N)}\prod_{n=1}^{\infty}
\Phi_N(q^n)^{24}$$
where  $\phi(N)$  is  Euler's totient  function,   $\Phi_N(X)$ is  the
$N^{th}$ cyclotomic polynomial and $q=e^{2\pi i z}$. This follows from
the M\"{o}bius inversion formula applied to the situation
$$\sum_{d|N} \log \Phi_d(X)=\log(1-X^N)$$
The inversion formula implies
$$\log(\Phi_N(X))=\sum_{d|N} \mu(d) \log(1-X^{\frac{N}{d}})$$
We also have 
$$\Phi_{N}(X)=\mathop{\prod_{k\;mod\;N}}_{(k,N)=1} (1-\xi^k X)$$
where $\xi=e^{\frac{2\pi i}{N}}$.

Define the {\bf $q$-logarithm} for $q=e^{2 \pi i z}$  as follows:
$${\bf log}_q(t)=\frac{1}{24}\log|qt|+\sum_{n=1}^{\infty} \log|1-q^nt|$$

Combining this with the formula for $L'(H^2(E \times E'),1)$, we get 
\begin{thm}[An `elliptic class number formula']
Let $ E,E',f, g$ be as before. We have 
$$L'(H^1(E) \otimes  H^1(E'),1)=\frac{-H(0)}{2} \mathop{\sum_{k\;mod\;N}}_{(k,N)=1}
\frac{1}{2 \pi i} \int_{X_0(N)} {\bf log}_q(\xi^k) f(q)\overline{g(q)}
\frac{dq}{q} \frac{d\bar{q}}{\bar{q}}$$

\end{thm}

\begin{proof}
From \ref{eis1} and \ref{lv1} we have 
\begin{equation}
L'_{f,g}(0))=-4\pi^2  i  \frac{1}{2  \pi i}  \int_{X_0(N)}
\phi(N)\log|q|\delta(f,g)
\end{equation}
$$-4  \pi^2 i \frac{1}{2  \pi i} \int_{X_0(N)}
\sum_{n=1}^{\infty}  \mathop{\sum_{k\;mod\;N}}_{(k,N)=1}  \log|1-\xi^k
q^n| \delta(f,g)$$
So it follows from the definition of ${\bf log}_q(t)$. Note the
similarity to \ref{cn1}.

\end{proof}

\section{Elements of $H_{\M}^{2n-1}(E_1 \times E_2
  \times ... \times E_n,\Q(n))_{\ZZ}$}

We can generalise this construction  to  prove the Beilinson and  Tate
conjectures for codimension $n-1$ cycles   on products of $n$  modular
elliptic  curves. We will  work  it  out  in detail  for $H^5_{\M}(E_1
\times E_2  \times E_3,\Q(3))$ of  the product  of $3$ elliptic curves
and remark how it generalises. It appears that all the cycles come from
exterior products.

Let $E_f$,$E_g$ and $E_h$ be $3$ modular elliptic curves corresponding
to the normalised  eigenforms $f$,$g$ and  $h$. Let $X=E_f \times  E_g
\times E_h$. From the description  of the real Deligne cohomology, one
can see that $H^5_{\D}(X_{/\R}.\R(3))$ is a 6 dimensional.

We  have to consider  the $L$-function of $H^4$ at  $s=2$. From the
Kunneth formula one can see that
$$L(H^4(X),s)=\zeta(s-2)^3L_{f,g}(s-2)L_{f,h}(s-2)L_{g,h}(s-2)$$
where $L_{f,g}(s)$ corresponds  to to the Rankin Selberg  convolution
of $f$ and $g$.

There are three cases that  we have to  consider, namely when all, two
or none of the three curves are isogenous.

\subsection{All the elliptic curves are isogenous}

In this case there are six elements of the Neron-Severi itself, namely
$$(x,e,e), (e,x,e), (e,e,x)$$
$$(x,x,e), (x,e,x), (e,x,x)$$
There is one more  `obvious' element,  namely  $(x,x,x)$, but there  is a
homology  relation between   these 7   elements,  giving  rise to  the
`modified diagonal cycle'. )

From the  calculations in     the  previous  section, one  has    that
$L_{f,g}(s-2), L_{g,h}(s-2)$ and $L_{f,h}(s-2)$ all have simple poles
are $s=3$ so that shows that the L function has a  pole of order $6$
as expected. From the  functional equation one sees that $L(H^4(X),s)$
is a nonzero rational number at $s=2$.

\subsection{Two are isogenous}

Assume only $E_f$ and $E_g$ are isogenous. In this case the conjecture
predicts that the rank of the motivic cohomology is $2$ while the rank
of the Neron-Severi is $4$. One has an exterior product map 
$$H^3_{\M}(E_f    \times    E_h,\Q(2))  \otimes    H^2_{\M}(E_g,\Q(1))
\longrightarrow H^5_{\M}(X,\Q(3))$$
From the earlier section, since $E_f$ and $E_h$ are not isogenous, one
has an element of $H^3(E_f  \times E_h,\Q(2))$ coming from the modular
parametrisation, and  one has the rational  point on $CH^1(E_g)$. This
gives an  element of $H^5_{\M}(X,\Q(3))$.   Similarly, using the other
pair of non-isogenous elements, one gets the other element.  These are
non-trivial as from the expression for the  $L$-function as a product,
the functions $L_{f,h}(s-1)$ and  $L_{g,h}(s-1)$ have simple zeroes at
$s=2$ and the  value of  $L''(H^4(X),s)_{s=2}$ is  the product of  the
regulators, which is the determinant of the regulator matrix.

\subsection{ None are isogenous }

In this case the conjectures  predict that there are three independent
elements of the  motivic cohomology and  three independent elements of
the Neron Severi.

We use the same argument as above, namely the exterior product and the
construction  in   the   previous section   to   construct  the  three
elements. 

\begin{rem}
  
  The same argument can also be  used to prove  this conjecture in the
  case of   products of {\em  modular}  curves as once  again  all the
  interesting elements come from $H^3_{\M}(X_1 \times X_2,\Q(2))$.

\end{rem}

\begin{tabular}{ll}
Ramesh Sreekantan & Srinath Baba \\
School of Mathematics             & Department of Mathematics   \\
Tata Institute of Fundamental Research       & McGill University  \\
Colaba  &  Montreal\\
Mumbai 400 005   & Quebec \\ 
India  & Canada  H3A 2K6 \\  
{\bf \verb+ramesh@math.duke.edu+}  & {\bf \verb+sbaba@math.mcgill.ca+} \\
         
\end{tabular}

\end{document}